\documentclass[11pt]{amsart}
\usepackage{latexsym,graphicx}

\numberwithin{equation}{section}
\theoremstyle{plain}

\theoremstyle{remark}

\theoremstyle{definition}

\newcommand{\D}{{\mathcal D}}
\newcommand{\E}{\mathcal E}

\newcommand{\G}{{\mathcal G}}

\newcommand{\K}{{\mathcal K}}
\renewcommand{\L}{{\mathcal L}}
\newcommand{\M}{{\mathcal M}}
\newcommand{\N}{\mathbb N}

\newcommand{\V}{{\mathcal V}}

\newcommand{\dist}{\operatorname{dist}}

\newcommand{\fp}{\operatorname{FP}}

\newcommand{\Int}{\operatorname{Int}}

\renewcommand{\span}{\operatorname{span}}

\newcommand{\supp}{\operatorname{Supp}}

\def\half{{1 \over 2}}
\def\la{\lambda}

\newcommand{\oa}{\overrightarrow}

\newcommand{\ol}{\overline}

\def\XXint#1#2#3{{\setbox0=\hbox{$#1{#2#3}{\int}$}
      \vcenter{\hbox{$#2#3$}}\kern-.5\wd0}}

\address{ Dept. of Mathematics, Rice  University, 6100 Main St., Houston, 77005 TX, U.S.A.
\\ {\sl E-mail address:}  {\bf harvey@rice.edu}}

\address{ Department of Mathematics, Stony Brook University, Stony Brook, NY 11790, U.S.A.
\\ {\sl E-mail address:}  {\bf blaine@math.sunysb.edu}}

\address{  Institute of Mathematics, Cracow University of Technology, Warszawska 24, 31-155
    Krak\'{o}w, Poland
\\ {\sl E-mail address:}  {\bf splis@pk.edu.pl}}

\begin{document}

\def\cal{\mathcal}

\font\tpt=cmr10 at 12 pt
\font\fpt=cmr10 at 14 pt

\font \fr = eufm10



\overfullrule=0in

\def\boxit#1{\hbox{\vrule
 \vtop{%
  \vbox{\hrule\kern 2pt %
     \hbox{\kern 2pt #1\kern 2pt}}%
   \kern 2pt \hrule }%
  \vrule}}

  \def\harr#1#2{\ \smash{\mathop{\hbox to .3in{\rightarrowfill}}\limits^{\scriptstyle#1}_{\scriptstyle#2}}\ }

 \def\GG{{{\bf G} \!\!\!\! {\rm l}}\ }

\def\GL{{\rm GL}}

\def\bll{I \!\! L}

\def\bra#1#2{\langle #1, #2\rangle}
\def\bbf{{\bf F}}
\def\bbj{{\bf J}}
\def\Jtn{{\bbj}^2_n}  \def\JtN{{\bbj}^2_N}  \def\JoN{{\bbj}^1_N}
\def\jt{j^2}
\def\jtx{\jt_x}
\def\Jt{J^2}
\def\Jtx{\Jt_x}
\def\bpp{{\bf P}^+}
\def\bpt{{\wt{\bf P}}}
\def\fsh{$F$-subharmonic }
\def\mo{monotonicity }
\def\jet{(r,p,A)}
\def\ss{\subset}
\def\sse{\subseteq}
\def\half{\hbox{${1\over 2}$}}
\def\smfrac#1#2{\hbox{${#1\over #2}$}}
\def\oa#1{\overrightarrow #1}
\def\dim{{\rm dim}}
\def\dist{{\rm dist}}
\def\codim{{\rm codim}}
\def\deg{{\rm deg}}
\def\rank{{\rm rank}}
\def\log{{\rm log}}
\def\Hess{{\rm Hess}}
\def\Hessyp{{\rm Hess}_{\rm SYP}}
\def\trace{{\rm trace}}
\def\tr{{\rm tr}}
\def\max{{\rm max}}
\def\min{{\rm min}}
\def\span{{\rm span\,}}
\def\Hom{{\rm Hom\,}}
\def\det{{\rm det}}
\def\End{{\rm End}}
\def\Sym{{\rm Sym}^2}
\def\diag{{\rm diag}}
\def\pt{{\rm pt}}
\def\Spec{{\rm Spec}}
\def\pr{{\rm pr}}
\def\Id{{\rm Id}}
\def\Grass{{\rm Grass}}
\def\Herm#1{{\rm Herm}_{#1}(V)}
\def\arr{\longrightarrow}
\def\supp{{\rm supp}}
\def\Link{{\rm Link}}
\def\Wind{{\rm Wind}}
\def\Div{{\rm Div}}
\def\vol{{\rm vol}}
\def\foral{\qquad {\rm for\ all\ \ }}
\def\fpsh{{\cal PSH}(X,\f)}
\def\Core{{\rm Core}}
\def\dis{f_M}
\def\Re{{\rm Re}}
\def\rn{\bbr^n}
\def\pp{\cp^+}
\def\plp{\cp_+}
\def\Int{{\rm Int}}
\def\cix{C^{\infty}(X)}
\def\Gr#1{G(#1,\rn)}
\def\Symn{{\Sym(\rn)}}
\def\SymN{{\Sym(\bbr^N)}}
\def\Gpn{G(p,\rn)}
\def\fd{{\rm free-dim}}
\def\SA{{\rm SA}}
 \def\cd{{\cal C}}
 \def\cdt{{\widetilde \cd}}
 \def\cm{{\cal M}}
 \def\cmt{{\widetilde \cm}}

\def\Theorem#1{\medskip\noindent {\bf THEOREM \bf #1.}}
\def\Prop#1{\medskip\noindent {\bf Proposition #1.}}
\def\Cor#1{\medskip\noindent {\bf Corollary #1.}}
\def\Lemma#1{\medskip\noindent {\bf Lemma #1.}}
\def\Remark#1{\medskip\noindent {\bf Remark #1.}}
\def\Note#1{\medskip\noindent {\bf Note #1.}}
\def\Def#1{\medskip\noindent {\bf Definition #1.}}
\def\Claim#1{\medskip\noindent {\bf Claim #1.}}
\def\Conj#1{\medskip\noindent {\bf Conjecture \bf    #1.}}
\def\Ex#1{\medskip\noindent {\bf Example \bf    #1.}}
\def\Qu#1{\medskip\noindent {\bf Question \bf    #1.}}
\def\Exercise#1{\medskip\noindent {\bf Exercise \bf    #1.}}

\def\HoQu#1{ {\AAA T\BBB HE\ \AAA H\BBB ODGE\ \AAA Q\BBB UESTION \bf    #1.}}

\def\pf{\medskip\noindent {\bf Proof.}\ }
\def\qed{\hfill  $\vrule width5pt height5pt depth0pt$}
\def\equdef{\buildrel {\rm def} \over  =}
\def\qedqed{\hfill  $\vrule width5pt height5pt depth0pt$ $\vrule width5pt height5pt depth0pt$}
\def\mathqed{  \vrule width5pt height5pt depth0pt}

\def\V{W}

\def\df{d^{\phi}}
\def\hk{\_{\rm l}\,}
\def\n{\nabla}
\def\w{\wedge}

\def\cu{{\cal U}}   \def\cc{{\cal C}}   \def\cb{{\cal B}}  \def\cz{{\cal Z}}
\def\cv{{\cal V}}   \def\cp{{\cal P}}   \def\ca{{\cal A}}
\def\cw{{\cal W}}   \def\co{{\cal O}}
\def\ce{{\cal E}}   \def\ck{{\cal K}}
\def\ch{{\cal H}}   \def\cm{{\cal M}}
\def\cs{{\cal S}}   \def\cn{{\cal N}}
\def\cd{{\cal D}}
\def\cl{{\cal L}}
\def\cp{{\cal P}}
\def\cf{{\cal F}}
\def\ccr{{\cal  R}}

\def\gerG{{\fr{\hbox{g}}}}
\def\gerB{{\fr{\hbox{B}}}}
\def\gerR{{\fr{\hbox{R}}}}
\def\p#1{{\bf P}^{#1}}
\def\vf{\varphi}

\def\wt{\widetilde}
\def\wh{\widehat}

\def\and{\qquad {\rm and} \qquad}
\def\arr{\longrightarrow}
\def\ol{\overline}
\def\bbr{{\mathbb R}}\def\bbh{{\mathbb H}}\def\bbo{{\mathbb O}}
\def\bbc{{\mathbb C}}
\def\bbq{{\mathbb Q}}
\def\bbz{{\mathbb Z}}
\def\bbp{{\mathbb P}}
\def\bbd{{\mathbb D}}

\def\a{\alpha}
\def\b{\beta}
\def\d{\delta}
\def\e{\epsilon}
\def\f{\phi}
\def\g{\gamma}
\def\k{\kappa}
\def\la{\lambda}
\def\o{\omega}

\def\s{\sigma}
\def\x{\xi}
\def\z{\zeta}

\def\D{\Delta}
\def\L{\Lambda}
\def\G{\Gamma}
\def\O{\Omega}

\def\bd{\partial}
\def\bdf{\partial_{\f}}
\def\lag{Lagrangian}
\def\psh{plurisubharmonic }
\def\ph{pluriharmonic }
\def\pph{partially pluriharmonic }
\def\omp{$\omega$-plurisubharmonic \ }
\def\ffl{$\f$-flat}
\def\PH#1{\widehat {#1}}
\def\lloc{L^1_{\rm loc}}
\def\dbar{\ol{\partial}}
\def\lp{\Lambda_+(\f)}
\def\lpp{\Lambda^+(\f)}
\def\bo{\partial \Omega}
\def\Ob{\overline{\O}}
\def\fc{$\phi$-convex }
\def\PSH{{ \rm PSH}}
\def\SH{{\rm SH}}
\def\totr{ $\phi$-free }
\def\BM{\lambda}
\def\Der{D}
\def\CH{{\cal H}}
\def\RH{\overline{\ch}^\f }
\def\pconv{$p$-convex}
\def\MA{MA}
\def\lagpsh{Lagrangian plurisubharmonic}
\def\hermsk{{\rm Herm}_{\rm skew}}
\def\PSHl{\PSH_{\rm Lag}}
 \def\ppsh{$\pp$-plurisubharmonic}
\def\fp{$\pp$-plurisubharmonic }
\def\fh{$\pp$-pluriharmonic }
\def\Symn{\Sym(\rn)}
 \def\ci{C^{\infty}}
\def\USC{{\rm USC}}
\def\fa{{\rm\ \  for\ all\ }}
\def\ppc{$\pp$-convex}
\def\cpt{\wt{\cp}}
\def\ft{\wt F}
\def\ob{\overline{\O}}
\def\Be{B_\e}
\def\K{{\rm K}}

\def\M{{\bf M}}
\def\N#1{C_{#1}}
\def\ds{Dirichlet set }
\def\dir{Dirichlet }
\def\Fa{{\oa F}}
\def\TR{{\cal T}}
 \def\LAG{{\rm LAG}}
 \def\ISO{{\rm ISO_p}}
 \def\Span{{\rm Span}}

\def\BB{1}
\def\CC{2}
\def\DD{3}
\def\EE{4}
\def\FF{5}

\def\CKNS{1}
\def\CTW{2}
\def\D{3}
\def\DF{4}
\def\GW{5}
\def\GN{6}
\def\PUP{7}
\def\DDI{8}
\def\DDR{9}
\def\GPSH{10}
\def\Survey{11}
\def\Bellman{12}
\def\Lagrangian{12}
\def\IDP{13}
\def\HLP{14}
\def\Plis{15}
\def\Pl{16}
\def\Pll{17}
\def\Plll{18}
\def\PlIll{18}
\def\Ri{19}
\def\TW{21}
\def\V{20}

\def\PSF{{\mathcal PSH}}

\def\E{E}
\def\bL{{\bf \Lambda}}
\font\headfont=cmr10 at 14 pt

\vskip .1in


\title[THE RICHBERG TECHNIQUE FOR SUBSOLUTIONS]
{THE RICHBERG TECHNIQUE FOR SUBSOLUTIONS}
 
\date{\today}
\author{ F. Reese Harvey,  H. Blaine Lawson, Jr.
and Szymon Pli\'s}
 \thanks{The second author was partially supported by the NSF and IHES,
 and the third author
was partially supported by the NCN grant 2013/08/A/ST1/00312.
}

\maketitle

\centerline{\sl Dedicated with great esteem to Karen Uhlenbeck.}

\vskip .3in
\centerline{\bf Abstract}
  \font\abstractfont=cmr10 at 10 pt
  
  {{\parindent= .3in\narrower \noindent
  
This note adapts the sophisticated Richberg technique for approximation in pluripotential
theory to the $F$-potential theory associated to a  general nonlinear convex subequation
$F \ss J^2(X)$ on a manifold $X$.  The main theorem is the following ``local to global'' result.
Suppose $u$ is a  continuous strictly $F$-subharmonic function  such that each point $x\in X$
has a fundamental neighborhood system consisting of domains for which a ``quasi'' form of  $C^\infty$
approximation holds.  Then for any positive  $h\in C(X)$  there exists a strictly $F$-subharmonic function
$w\in C^\infty(X)$ with $u< w< u+h$.  Applications include all convex constant coefficient subequations on $\rn$,
 various nonlinear subequations on complex and almost complex manifolds, and many more.

}}

\vskip .6in
\centerline{\bf Table of Contents}

\medskip

\hskip1in
 \BB.     Introduction.

\hskip1in
 \CC.    The Proof of the  Local to Global Theorem.

\hskip1in
 \DD.    Verifying the Local Hypothesis.

\hskip1in
 \EE.    Smoothing the Maximum Function.

\hskip1in
 \FF.    An Example -- Non-convex Subequations.


\vfill\eject


\medskip
\noindent{\headfont \BB.\  Introduction.}

The point of this paper is to extend the classical Richberg technique 
in pluripotential theory (cf. [\Ri] or [\D, Lemma 5.17/Cor. 5.19/Thm. 5.21]) to 
subsolutions of any
convex subequation. Instances of this have appeared in [\GW],  [\PUP],  [\GPSH],  [\Pl], [\Pll], [\Plll], and [\V].
Here we present a quite broad result.

The main idea is very general and can be formulated in the context of the
potential theory associated with any convex subequation 
$F\ss J^2(X)$\footnote{$F$ is a subequation if it satisfies  positivity $F+\cp\ss F$, negativity $F+\cn\ss F$,
and a mild topological condition (T) [\DDR, p. 416]. Here $\cp_x$ is the set of 2-jets of non-negative functions with critical value xero at $x$, and $\cn_x$ are the jets of non-positive constant functions.
} 
on a manifold $X$.  The space $F(X)$ of {\bf subsolutions} or $F$-{\bf subharmonic functions}
consists of all upper semicontinuous functions $u:X\to \bbr\cup \{-\infty\}$
with the property that for any $x\in X$ and any test function\footnote{A $C^2$-function $\vf$
such that $u-\vf$ has a local max of zero at $x$}
 $\vf$ for $u$ at $x$, the 2-jet  at $x$ of $\vf$ lies in $F$ (denoted $J^2_x \vf\in F$).

A smooth function $v\in F(X)$ is said to be {\bf $F$-strict}   if $J_x^2 v \in\Int F$ for all $x\in X$.
A general $u\in F(X)$  is said to be $F$-strict if, given $\vf\in C^\infty_{\rm cpt}(X)$, the
function $u+\e \vf$ is $F$-subharmonic for all $\e>0$ sufficiently small.

It is useful to further refine this notion as follows.  Given a strictly positive  $g\in C(X)$,
a continuous function $u\in F(X)$ is said to be {\bf $g$-strict} if $u$ is $F^g$-subharmonic   where
$F^g$ is the subequation with fibres
$$
F_x^g \ =\ \{ J \in J^2_x(X) : \dist(J, \sim F_x)\geq g(x)\}.
$$
One has that
$$
\text{$u$ is $F$-strict $\qquad\iff\qquad$  $u$ is $g$-strict for some $g>0$ on $X$. }
\eqno{(\BB.1)}
$$
See Corollary 7.6 in [\DDR] for the proof of the implication $\Leftarrow$.
(We note that the above definition of {\sl $F$-strict} is implicit in [\DDR].)
The proof of the  implication $\Rightarrow$ is left to the reader.

The appropriate local approximation hypothesis is the following  form of 
``quasi''-approximation. Let $F^{\rm strict}_{\rm cont}(X)$ denote the continuous $F$-strict functions on $X$.

\noindent
{\bf Definition \BB.1.}  Suppose $u\in F^{\rm strict}_{\rm cont}(X)$.

\noindent
(a)\ \   Given a domain
$\O\ss\ss X$, we say that {\bf quasi $C^\infty$ approximation holds for $u$ on $\O$} if 
for all compact sets $K\ss \O$, there exists 
$v\in C(\ob) \cap  F^{\rm strict}(\O) \cap \ci(\O)$
such that:
\medskip

\centerline{(A) \ $u<v$ on $K$ \and  (B) \ $u>v$ on $\bo$.}

\medskip
\noindent
(b)\ \ If for each $h\in C(X)$, with $h>0$,  
 there exists $$w\in C^\infty(X)\cap F^{\rm strict}(X)$$ which satisfies:
$$
u\ <\ w\ <\ u+h 
$$
on $X$, then we say that {\bf Richberg approximation holds} for $u$ on $X$.

Now we can adapt the Richberg technique to subequations to prove the following 
local-to-global result. (We shall always assume that the subequation $F$ is convex, 
i.e., the fibre at each point is convex.)

\Theorem{\BB.2} {\sl
Suppose $u\in  F^{\rm strict}_{\rm cont}(X)$ and that each point $x\in X$ has a fundamental neighborhood
system consisting of domains $\O$ for which quasi  $C^\infty$ approximation holds for $u$.
Then  Richberg approximation holds for $u$ on $X$.}

\noindent
{\bf Remark \BB.3.}  If the local  $C^\infty$ approximators $v_\a$ to $u$ used in the proof can be
chosen to be $g$-strict, for positive $g\in C(X)$, then the 
function $w$ constructed in the proof can be chosen to be $g'$ strict for any $0<g'<g$.  (For the proof
see Lemma \CC.2 below.)

We now discuss the  applications of Theorem \BB.2.  The proofs, together with some further
results, are given in Section \DD.

We begin with a  subequation on $\rn$ which is both
 constant coefficient  and convex.

\Theorem{\BB.4}
{\sl
Suppose $F$ is a  
constant coefficient, convex  subequation on $\rn$ and $u \in F^{\rm strict}_{\rm cont}(X)$
for some open subset $X\ss \rn$. Then   Richberg approximation holds for $u$ on $X$.}

This result extends the pure second-order case provided by Thm.\ 9.10 in  [\GPSH].

It applies to all convex subequations 
$$F\  \ss\  J^2(\rn)\  =\  \rn\times(\bbr\oplus\rn\oplus\Symn).$$
Among these are many highly degenerate subequations which 
are geometrically interesting -- for example, the 
subequations $F(\GG)$   defined by a closed subset $\GG\ss G(p, \rn)$
of the Grassmannian of $p$-planes in $\rn$.  These include the subequations coming from calibrations,
and those associated to Lagrangian planes (or more generally isotropic planes) in $\bbc^n$
(see  [\GPSH], [\Lagrangian] for more details). 



For reduced convex subequations $F$ on a  manifold $X$ we assume   a mild  condition
on coordinate balls  (see Definition \DD.12 and Lemma \DD.13), and obtain the following two results.

The first is for $F$ a cone, and is based on solving the $C^\infty$ homogeneous Dirichlet problem.
In the second result, the cone hypothesis is dropped, but sufficient ``monotonicity'' is assumed.
This second case is again  based on solving the $C^\infty$  Dirichlet problem, but this time for the
inhomogeneous equation.

\Theorem{\BB.5. (Convex Cone Subequations)}
{\sl   Let  $F$ be a  reduced\footnote{independent of the value of the function} 
convex cone  subequation on a manifold $X$.
Suppose  that $F$ satisfies the coordinate ball condition in Definition \DD.12.
 Assume that    
on all sufficiently small coordinate balls  $\O$ for a covering family of local coordinates on $X$
one has that  the  $C^\infty$ homogeneous  Dirichlet problem for $F$ is uniquely 
 solvable on $\O$ (see Def.\ \DD.6).
 Then for all $u \in F^{\rm strict}_{\rm cont}(X)$
  Richberg approximation holds for $u$ on $X$.
  }

\Theorem{\BB.6. (Subequations with a Good Monotonicity Cone)}
{\sl   Let  $F$ be any  reduced  subequation with a monotonicity cone $M$ on a manifold $X$.
  That is, $M$ is a reduced convex cone 
subequation such that the fibre-wise sum satisfies
$$
F+M \ \ss\ F.
$$
Now suppose  that $M$  satisfies the coordinate ball condition in Definition \DD.12.
 Assume that    
on all sufficiently small coordinate balls  $\O$ for a covering family of local coordinates on $X$
one has that  the  $C^\infty$ inhomogeneous  Dirichlet problem  for $M$ is uniquely 
 solvable on $\O$ (see Def.\ \DD.8).
 Then for all $u \in F^{\rm strict}_{\rm cont}(X)$
  Richberg approximation holds for $u$ on $X$.
  }

A special case of this theorem is when $F$ is a convex cone subequation and one takes $M=F$.
See Remark \DD.14.

Note that Theorem \BB.5 is not a special case of Theorem \BB.6 since in \BB.6 $\psi\equiv 0$ is excluded.

\Ex{\BB.7. (Plurisubharmonics on an Almost Complex Manifold)}  Let $(X,J)$ be an almost complex $n$-manifold, and $\cp(J)$ the subequation
defining the $J$-plurisubharmonic functions by the condition $J^2_x \vf \in F_x(J)$ iff $i\partial\dbar \vf \geq0$.
Fixing a volume form $\b$ on $X$, there is a natural operator $f(J^2\vf) = (i\partial\dbar \vf)^n/\b$.
It was shown in [\Plis], [\CTW] that the $\ci$  inhomogeneous Dirichlet problem is uniquely solvable
on small balls in local coordinates. This says that for a smooth function $\psi >0$, the subequation
$\bbf \equiv \{J: f(J)\geq \psi\}$ satisfies Definition \DD.6.   Richberg approximation follows from Theorem \BB.5.
This result  was first established
in [\Pll].

\Ex{\BB.8. (Complex Hessian Equations)}  Suppose $(X,\o)$ is a K\"ahler $n$-manifold, 
and $F$ is the complex $m$-Hessian subequation:
$J^2_x \vf\in F_x$ iff $(i\partial\dbar \vf)^k \wedge \o^{n-m}\geq0$,  for all    $1\leq k\leq m$. 
Fixing a volume form $\b$ as above
gives an operator $f(J^2 \vf) \equiv (i\partial\dbar \vf)^m\wedge \o^{n-m}/\b$
for which the $\ci$ inhomogeneous Dirichlet Problem is uniquely solvable on small
coordinate balls (see [\Plll] or [\GN]). Hence Theorem \BB.6 applies to yield   Richberg approximation on $X$.

\Ex{\BB.9.  (Work of Verbitsky and Greene-Wu)} Suppose again that  $(X,\o)$ is a K\"ahler $n$-manifold, and let $F$ be the subequation:
$J^2_x \vf\in F_x$ iff $i\partial\dbar \vf \wedge \o^{k}\geq0$ for fixed $k$, $0\leq k\leq n-1$.
Equivalently this is the subequation given by demanding that the trace of the complex hessian is non-negative
on every tangent complex  $(k+1)$-plane.
The Richberg Approximation Theorem was proved in these cases by M. Verbitsky [\V].
His argument was based in part on a general local-to-global result
due to R.\ E.\ Greene and H.\ Wu [\GW].

In Theorems \BB.5 and \BB.6 the subequation is required to be convex (and reduced).
An example of a non-convex subequation where nevertheless   Richberg approximation holds is discussed
in Section \FF.

\Remark{\BB.10. (Justifying the Definitions)}
We point out  that in Definitions \DD.6. \DD.8   the $\ci$-regularity 
for solutions to the Dirichlet problem  is only assumed to hold
on $\O$ (given a $\ci(\bo)$ boundary function).  Consequently, Theorem \BB.5
applies to many more of the examples in the literature than it would if we required
$\ci$-regularity on $\ob$.

On the other hand, the weakest hypothesis required for adapting the local-to-global 
technique to subequations (as provided by Theorem \BB.2)  is the notion of {\sl quasi $\ci$ approximation}
given in Definition \BB.1(a).   Theorem \BB.5 then follows from Theorem \BB.2 by 
proving that regularity for the Dirichlet problem implies quasi $\ci$ approximation.

We  introduce another weak notion -- that of {\sl approximate $\ci$-regularity} (Definition \DD.4) --
which is implied by the full $\ci$-regularity for the homogeneous or inhomogeneous Dirichlet problem
(Propositions \DD.7, \DD.9), but in principle is much easier to establish.  This approximate regularity
implies our weakest hypothesis  of quasi $\ci$ approximation (Lemma \DD.2).

For the constant coefficient subequations in Theorem \BB.4, standard $\ci$ approximation holds
by the usual convolution techniques (Lemma \DD.3). This reduces the proof of Theorem \BB.4 to 
Theorem \BB.2 by showing that $\ci$ approximation implies quasi $\ci$ approximation (Lemma \DD.2).


\medskip
\noindent{\headfont \CC.\  The Proof of the Local to Global Theorem.}

The proof of Theorem \BB.2 relies on a regularization of the maximum function $M(t) \equiv \max\{t_1,...,t_m\}$
on $\bbr^m$ by convolution based on an approximate identity which is a product.
More specifically, choose $\vf(s)\in C^\infty_{\rm cpt}(\bbr)$ with $\vf\geq0$, $\supp \vf\ss[-1,1]$,
and $\int_{\bbr} \vf(s) \,ds =1$.  For each $\e=(\e_1,...,\e_m)$ with $\e_j>0$, 
set $\vf(y)\equiv \vf(y_1)\cdots \vf(y_m)$ and 
$$
\vf_\e(t) \ \equiv\ \vf\left( {t_1\over \e_1}  \right) \,\cdots \, \vf\left( {t_m\over \e_m}  \right) {1\over \e_1\cdots \e_m}
$$
Define
$$
M_\e(t)  \ \equiv\ \int_{\bbr^m} M(t+y) \vf_\e(y)\,dy
\ =\  \int_{\bbr^m} M(t+\e y) \vf(y)\,dy,
\eqno{(\CC.1)}
$$
using an ``abuse of notation''  $\e y = (\e_1  y_1,...,\e_m y_m)$.  Assume also that $\vf$ is an even function, so that
$$
 \int_{\bbr^m} y_j \, \vf_\e(y) {dy_1 \over \e_1} \,\cdots\,  {dy_m \over \e_m}  \ =\ 0.
\eqno{(\CC.2)}
$$
This ensures that the convolution of  $\vf_\e$ with each $t_j$ equals zero.

\noindent
{\bf Properties \CC.1.} {\sl

(1) $M_\e(t)$ is increasing in all the variables, smooth and convex on $\bbr^m$, 
and invariant  under permutations of the variables.

(2) $M_\e(t+se) = M_\e(t)+s$ where $e=(1,...,1)$, and hence 
$$\sum_{j=1}^m {\partial M_\e  \over \partial t_j} \equiv 1.$$

(3) $M(t) \leq  M_\e(t) \leq M(t+\e)$.

(4) If $t_j+\e_j \leq \max_{i\neq j} \{t_i-\e_i\}$, then 
$$
M_\e(t)  \ = \ M_{(\e_1,...,\widehat{\e_j},...,\e_m)} (t_1, ... , \wh{t_j}, ... ,t_m).
$$
}

The proof of these properties will be discussed 
in Section \EE, along with a proof of the next result.

\Lemma{\CC.2} {\sl
Suppose that $F$ is a convex subequation (not necessarily reduced nor a cone) 
on a manifold $X$.
Suppose $g\in C(X)$, $g>0$.  If $u_1, ... , u_m \in F^g(Y)\cap C^\infty(Y)$
for  $Y\ss X$, then we have that}
$$
M_\e(u_1,...,u_m) \in F^{g-M(\e)}(Y)\cap C^\infty(Y).
$$

\noindent
{\bf Proof of Theorem \BB.2.} Pick a locally finite open cover $\{\O_\a'\}$  of $X$
consisting of precompact domains in $X$ with the property that quasi $C^\infty$ approximation
 holds for $u$ on each $\O_\a'$.
By a standard result in topology we can choose a subordinate open covering  $\{V_\a\}$ of $X$
with  $K_\a \equiv \overline{V_\a} \ss\O_{\a}'$ for all $\a$.
By the quasi $\ci$-approximation property for $u$ on $\O_\a'$, there exist
functions $v_\a \in C(\overline{\O}_{\a}') \cap \ci(\O_{\a}') \cap F^{\rm strict}(\O_{\a}')$
such that 
$$
(A)\ \ u <v_\a \quad {\rm on} \ \ K_\a
\and
(B)\ \  u >  v_\a \quad {\rm on} \ \ \partial  \O_\a'.
$$
Now choose open sets $\O_\a \supset K_\a$ with $\overline{\O}_\a\ss\O_{\a}'$.
Then 
$$
v_\a \in \ci({\rm nb\, } \overline{\O}_\a) \cap F^{\rm strict} ({\rm nb\, } \overline{\O}_\a)
$$
(where ${\rm nb\, } \overline{\O}_\a$ denotes a neighborhood of $\overline{\O}_\a$).
Furthermore, by choosing each $\O_\a$ sufficiently large,   we can arrange that 
$$
(A)\ \  u <v_\a  \quad {\rm on} \ \ K_\a
\and
(B)\ \  u >  v_\a \quad {\rm on} \ \ \partial  \O_\a.
$$

\noindent
In addition, by choosing the $\O_\a'$ sufficiently small we can assume that
\medskip

\centerline{(C) \ $\sup_{\ob_\a} u < \inf_{\ob_\a} (u+h)$.}

Next we choose $\e_\a$ so that 
\medskip

\centerline{(A$'$) \ $u<v_\a-\e_\a$ on $K_\a$ \and  (B$'$) \ $u>v_\a + \e_\a$ on $\bo_\a$.}

We can define the global function $w$ on $X$ (using the fact that $M_\e$ is invariant
under permutations) by
$$
w(x) \ \equiv\ M_\e\{v_\a(x) : x\in \O_\a\}
\eqno{(\CC.3)}
$$

It remains to prove that $w$ satisfies all the properties in Theorem \BB.2.
Note that the smoothness (even the continuity) of $w$ is completely
unclear from this definition. 

We now fix a smooth function $g'$ with $0<g'<g$, and we reduce the $\e_\a$ 
so that $0< \e_\a \leq \sup_{\O_\a}(g-g')$.  We shall first  give a proof that:
\medskip

\centerline{ $w\in C^\infty(X) \cap F^{g'}(X)$.}

\noindent
 Given $x\in X$ there exists a neighborhood of $x$  where
$$
M_\e\{v_\a(y) : y\in \O_\a\}\ =\ M_\e\{v_\a(y) : y\in \O_\a \ {\rm and}\  x \notin \bo_\a\}
\eqno{(\CC.4)}
$$
and hence
$$
w(y) \ =\ M_\e\{ v_\a(y) : x\in \O_\a\} \ \ \text{in a neighborhood $Y$ of $x$.}
\eqno{(\CC.5)}
$$
To prove (\CC.4) note that if $x\in\bo_\a$, then since $\{K_\b\}$ covers $X$, there exists 
$\b$ with $x\in K_\b \ss\O_\b$, which by (B) and (A) implies that
$$
v_\a(x) \ <\ u(x) \ <\ v_\b(x).
\eqno{(\CC.6)}
$$
Thus $v_\a < v_\b$ is a neighborhood of $x$.  Now by (\CC.5)
we see that $w$ is smooth on the neighborhood $Y$.
Furthermore,  Lemma \CC.2 applies to prove the $g'$-strictness of $w$
discussed in Remark \BB.3.

Next we prove that $u<w$.
Since each point  $x$ is contained in $K_\a \ss\O_\a$ for some $\a$, condition (A) implies
that $u(x) < v_\a(x) \leq M\{v_\b(x) : x\in \O_\b\}$, which by the first part of Property (3) is 
$\leq M_\e\{v_\b(x) : x\in \O_\b\} \equiv w(x)$.

Finally we prove that $w<u+h$.
By (B$'$) and the hypothesis (C) we have
$$
\sup_{\bo_\b} (v_\b+\e_\b) \ \leq\ \sup_{\bo_\b} u \ \leq \ \inf_{\ob_\b} (u+h).
\eqno{(\CC.7)}
$$
Now 
$$
w(x) \ \equiv \  M_\e\{v_\b(x)  : x\in \O_\b\} 
\ \leq  \   M\{v_\b(x) + \e_\b : x\in \O_\b\} 
$$
by the second part of  Property (3).  Since each $v_\b +\e_\b$ satisfies the Maximum Principle on $\O_\b$,
we have 
$$
w(x) \ \leq \ M\{\sup_{\bo_\b} (v_\b+\e_\b) : x\in \O_\b\}.
$$
Finally, by (\CC.7) we have $w(x) \leq u(x) + h(x)$. \qed

\vskip.3in

\noindent 
{\headfont \DD. Verifying the Local Hypothesis.}

In this section we give the proofs of Theorems \BB.4, \BB.5 and \BB.6
by verifying the local hypothesis of quasi $\ci$-approximation
so that the local-to-global  Theorem \BB.2 can be applied.  Two distinct methods are used.
  The first will apply to constant coefficient equations in $\rn$, yielding Theorem \BB.4.
The second will apply to a number of important nonlinear equations on manifolds where
the Dirichlet problem is sufficiently well understood, and will yield Theorems \BB.5 and \BB.6.

\medskip
\noindent
{\bf METHOD 1.}  Not surprisingly our first method is based on standard convolution.  
It also employs  the following notion.

\noindent
{\bf Definition \DD.1}
Given $u\in F^{\rm strict}_{\rm cont}(X)$ and a domain $\O\ss\ss X$, we say that 
{\bf $C^\infty$ approximation holds for $u$ on $\O$} if there exists a constant $c>0$ and a sequence
$\{v_j\}$ of functions $v_j$ which are $\ci$ and $F^c$-strict on a neighborhood of $\ob$ and converge
uniformly to $u$ on $\ob$.

Quite naturally,  $C^\infty$ approximation  implies quasi $C^\infty$ approximation , but there is
something to check.

\noindent
{\bf Lemma \DD.2} {\sl
Suppose  $C^\infty$ approximation  holds for  $u\in F^{\rm strict}_{\rm cont}(X)$  on a domain $\O$.  Then
quasi $C^\infty$ approximation  holds for  $u$ on $\O$.
}

\noindent
{\bf Proof.} 
Suppose $K$ is a compact subset of $\O$.  Choose a defining function $\rho\in C^\infty(\ob)$
for $\bo$ (neither the $F$-subharmonicity of $\rho$ nor the smoothness of $\bo$ are required in this lemma).
Then choose $s>0$ such that $\rho+s <0$ on $K$.  Consider the function
$$
{\overline v}_j \ \equiv \ v_j -\d (\rho+s) \quad \text{on a neighborhood of } \ \ob.
\eqno{(\DD.1)}
$$
Since each $v_j$ is $c$-strict, 
 given $0< c' <c$, there exists $\d_0>0$ such that 
 for $\d \leq \d_0$
 $$
 {\overline v}_j \ \ \text{is $c'$-strict on $\ob$ for all $j$}.
 \eqno{(\DD.2)}
$$
Now ${\overline v}_j - u = -\d(\rho+s) + v_j -u$.  Since $-\d(\rho+s)>0$ on $K$ and 
$-\d(\rho+s)<0$ on $\bo$, by choosing $v_j-u$ small enough on $\ob$ both (A) and (B) hold for 
${\overline v}_j$.\qed

This result can be used in the following case.

\vskip.2in

\centerline{\bf Convex Constant Coefficient Subequations.}

\Lemma{\DD.3}  {\sl
Suppose $F$ is a convex constant coefficient subequation on $\rn$.
Given $u\in F^{\rm strict}_{\rm cont}(X)$ and a domain $\O\ss\ss X$, 
$C^\infty$ approximation holds for $u$ on $\O$.
}

Standard convolution can be used, but there are things to prove, and a basic result
of [\Bellman] is necessary.

\noindent
{\bf Proof.}  There exists a constant $c>0$ such that $u$ is $c$-strict on a neighborhood of 
$\ob$ (see Section 7 in [\DDR]), i.e., $u$ is $F^c$-subharmonic on a neighborhood of $\ob$.
For any convex set in  euclidean space the set of points of distance $\geq c$ to the complement
is also convex since minus the log of the distance is convex.  Thus $F^c\ss F \ss \bbj^2$ is convex.
Also note that $F^c = \overline{\Int F^c}$
  The condition (P)   holds for $F^c$ by Lemma 7.3 in [\DDR].
In summary, $F^c$ is a convex constant coefficient subequation.  Consequently, 
because of the results proved in [\Bellman],
there is an equivalent
approach to subsolutions using distribution theory.  Now using the distributional approach, 
one shows that   $C^\infty$ approximation holds on a neighborhood of $\ob$ via standard convolution. \qed

\noindent    
{\bf Proof of Theorem \BB.4.} Combining Lemmas \DD.3 and \DD.2  proves that quasi $\ci$ approximation holds,
so that Theorem  \BB.2 applies. \qed

\medskip

\medskip
\noindent
{\bf METHOD 2.} 
Our second method is based on assuming certain regularity hypotheses concerning the Dirichlet Problem.
Since the discussion is local  we can   assume that the subequation $F$ is defined  on an open set
$X\ss\rn$.  
 We fix a domain $\O\ss\ss X$ with smooth boundary $\bo$ and let $\rho \in \ci(\ob)$
denote a defining function for $\bo$.  We assume that 
existence and comparison hold for the Dirichlet Problem (DP) on $\O$, and then
consider  the following form of ``regularity''.

\noindent
{\bf Definition \DD.4}  We say that   {\bf  approximate $\ci$-regularity for the (DP) holds on $\O$}
if, for each boundary function $\vf\in\ci(\bo)$, the  solution $H$ to the (DP)
can be uniformly approximated on $\ob$ by a sequence of functions 
$\{v_k\}\ss C(\ob)\cap F^{\rm strict} (\O) \cap \ci(\O)$.

Note that $v_k$ is not required to be $\ci$ on $\ob$.  This kind of regularity is enough for our purposes.
Our next result says that the quasi $\ci$-approximation hypothesis in Theorem \BB.2 can be replaced 
by this approximate $\ci$-regularity hypothesis as long as the subequation is reduced (Thm. \BB.2$'$ below).

\Lemma{\DD.5}
{\sl  
Suppose that $F$ is a reduced subequation on $X$.
If approximate $\ci$-regularity for the (DP) holds on $\O$, then quasi $\ci$-approximation
holds for all $u\in F^{\rm strict}_{\rm cont}(X)$  on $\O$. 
}

\pf   Suppose a compact subset $K\ss \O$ is given.  
Since $u$ is $F$-strict on $X$, for 
 $\e>0$ sufficiently small   
$$
u- \e \rho \quad  \text{ is $F$-subharmonic on a neighborhood of }\ \ob.
\eqno{(\DD.3)}
$$
Now choose a constant $s$ with  $0< s < \inf_K (-\rho)$ so that  
$$
\rho + s \ <\ 0\quad  \text{ on}\ K.
\eqno{(\DD.4)}
$$
Finally,  choose $\vf\in C^\infty(\bo)$ with
$$
u- \e s\ <\ \vf\ <\ u\quad  \text{ on}\ \bo.
\eqno{(\DD.5)}
$$
Let $H$ be the solution to the (DP) on $\O$ with boundary values $\vf$.  First
we want to prove that

\hskip .5in (A) \ \ $u\ < \  H$ \ \ on $K$
\and
 (B) \ \ $H\ < \ u$\ \  on $\bo$

\noindent
Assertion (B) follows since $\vf < u$ on $\bo$  by (\DD.5).
Now set $w\equiv u - \e (\rho+s)$ and note that $u<w$ on $K$ by  (\DD.4).
To finish the proof of  (A) we show that $w\leq H$ on $\ob$ (and therefore on $K$)
by showing that $w$ is in the Perron family $\cf(\vf)$.  
Note that on $\bo$ we have $w= u-\e s$, which is $<\vf$ by (\DD.5).
By (\DD.3) and the hypothesis that $F$ is reduced, $w$ is $F$-subharmonic on $\O$.  Hence $w\in \cf(\vf)$ as claimed.

Now  $H$ is neither strict nor $\ci$, however by the hypothesis there
exists a sequence $\{v_k\}\ss C(\ob)\cap F^{\rm strict} (\O) \cap \ci(\O)$ converging to $H$
uniformly on $\ob$.  Hence (A) and (B) hold with $H$ replaced by $v_k$ for large $k$, 
 proving that quasi $\ci$-approximation holds for $u$ on $\O$.\qed

\Theorem{\BB.2$'$}  {\sl
Suppose $F$ is a reduced convex cone subequation on a manifold $X$.  
Suppose $u \in F^{\rm strict}_{\rm cont} (X)$ and that each point $x\in X$ has a fundamental 
neighborhood system consisting of domains for which approximate $C^\infty$ regularity
holds for the (DP).  Then Richberg approximation holds for $u$ on $X$.}

We now discuss two cases where approximate $\ci$-regularity holds
so that Theorem \BB.2$'$ can be applied.

\bigskip
\centerline{\bf $\ci$ Regularity for the Homogeneous $F$-Dirichlet Problem.}

{\bf Definition {\DD.6}.} We say that the  {\bf    $C^\infty$ homogeneous  $F$-Dirichlet problem  (DP)  is uniquely 
 solvable on $\O$}  if for each  $\vf\in \ci(\bo)$, there exists a unique $h\in C(\ob)$
 such that

(a) \ \ $h$ is $F$-harmonic on $\O$,

(b) \ \ $h\bigr|_{\bo} = \vf$,

(c) \ \ $h \in \ci(\O)$ (but not necessarily in $\ci(\ob)$).

In the next result we assume that $F$ is a convex cone subequation.

\noindent
{\bf  Proposition {\DD.7}.} {\sl
Assume that there exists  a defining function $\rho$ for $\bo$
 which  is strictly $F$-subharmonic on a neighborhood of $\ob$.
If the $\ci$ homogeneous  Dirichlet problem is uniquely 
 solvable on $\O$,  then approximate $\ci$-regularity for the (DP) holds on $\O$.
}

\noindent
{\bf Proof.}
Suppose $h$ is the solution of the (DP) on $\O$ with boundary values $\vf \in \ci(\bo)$.
Then the sequence $v_k  \equiv h +{1\over k} \rho$ 
 uniformly approximates $h$ on $\ob$ and  $v_k\in  C(\ob)\cap F^{\rm strict} (\O) \cap \ci(\O)$ 
for each $k$.\qed


\bigskip
\centerline{\bf $\ci$ Regularity for the Inhomogeneous $F$-Dirichlet Problem}
\centerline{\bf Given a Monotonicity Cone.}

We shall now widen the  interesting case above to a broader family of subequations
by dropping the assumption that $F$ is a cone.  However, we
 now assume  that the subequation $F$ admits a {\bf monotonicity cone} $M$, 
that  is, $M$ is a convex cone subequation on $X$ for which the fibre-wise sum satisfies
$$
F+M \ \ss\ F.
$$

In addition to the subequations $F$ and $M$ we require weakly elliptic operators
\footnote{\ On a manifold $X$ weak ellipticity means  $f(J+J')\geq f(J)$ for $J\in F$ and $J'\in \cp$ 
(see footnote 1).} 
$f\in C(F)\cap C^\infty(\Int F)$ and $g\in C^\infty(M)$ which are compatible with the 
subequations, by which we mean  that
$$
\text{
$f\geq 0$ on $F$ and  $\partial F = \{f=0\}$ \quad and \quad 
$g \geq 0$ on $M$ and  $\partial M= \{g=0\}$.
}
$$
Finally, we assume that 

(1)  \ \ There is a constant $\d>0$ such that  if $J\in F$ and $J'\in M$, then $$f(J+J') \geq    \d g(J').$$

Next we show that interior $C^\infty$ regularity for the inhomogeneous
$F$-Dirichlet Problem implies that  approximate $\ci$-regularity for the $F$-(DP) holds on $\O$.

\noindent
{\bf Definition \DD.8.}
Assume that there exists a defining  function $\rho\in \ci(\ob)$ for $\bo$,
which is strictly $M$-subharmonic on $\O$.
We say that the {\bf $C^\infty$ inhomogeneous Dirichlet Problem (IHDP)  is uniquely solvable  on $\O$} 
if for all $\vf\in C^\infty(\bo)$ and $\psi \in C^\infty(\ob), \psi > 0$,
there exists a unique $v\in C(\ob)$ satisfying

(a) \ \ $ f(J^2 v)= \psi \quad $ on $\O$,

(b) \ \ $v\bigr|_{\bo} =\vf$,

(c) \ \ $v\in C^\infty(\O)$ (not necessarily in $C^\infty(\ob) $), and

(d) \ \  $v$ is the Perron function, i.e., 
$$
v(x) \ =\  \sup\{ u(x) : u\in\USC(\ob) \cap F(\O), \ u\bigr|_{\bo}\leq \vf, \
 {\rm and} \ f(J^2_{\rm red} u) \geq \psi \ {\rm on}\ \O\}
$$

\noindent
{\bf Proposition {\DD.9}.}
{\sl
Assume that there exists  a defining function $\rho$ for $\bo$
 which  is strictly $M$-subharmonic on a neighborhood of $\ob$.
If the $C^\infty$  (IHDP) is uniquely solvable on  $\O$, then approximate $C^\infty$-regularity for the (DP)
 holds on $\O$.
}

\pf
Suppose $H$ is the solution to the homogeneous $F$-(DP) with boundary values $\vf\in \ci(\bo)$.
Let $v_\e$ be the solution to the (IHDP) $f(J^2 v_\e) = g(J^2(\e \rho))$ with boundary values $\vf$.
Since $g(J^2(\e \rho))>0$ on $\O$ (by the strict $M$-subharmonicity of $\rho$), we see that  $v_\e$ is $F$-subharmonic.  Since $H$  is the Perron function for the (DP) and $v_\e =\vf$ on $\bo$, we have
$$
v_\e \leq H.
$$
  On the other hand, since
$f$ satisfies (1), 
$$
f\left(J^2(H+\e \rho)\right)\ \geq\  \d g(J^2(\e \rho)).
$$
Therefore $H+\e\rho$ is in  the Perron family for the (IHDP) and equals $\vf$ on $\bo$, so that
$$
H+\e\rho \leq v_\e.
$$
Thus, $0\leq H-v_\e \leq -\e\rho$ and so $v_\e \to H$ uniformly on $\ob$.
By our hypothesis, each $v_\e$ is $\ci$ on $\O$.\qed

\noindent
{\bf Note \DD.10.}  The solution $H$ to the homogeneous Dirichlet problem in the
proof of Proposition \DD.9, can simply be a continuous viscosity solution.

\Ex{\DD.11} Suppose $(X, \o)$ is a K\"ahler manifold.  Let $F$ be defined by the condition
$\{u : i\partial \dbar u +\o\geq0\}$,  and  define $f$ to be the complex determinant of
$i\partial \dbar u +\o$.  Let  $M$ be given by  $\{u : i\partial \dbar u\geq0\}$, and set
$g$ = the complex determinant of  $i\partial \dbar u$.

\medskip

To complete the proofs of Theorems \BB.5 and \BB.6 we need a final lemma.
The hypothesis of approximate $\ci$ regularity in Theorem \BB.2$'$ only needs to hold
on a family of small domains about each point.  For that reason we shall make the following 
definition below.

We first recall  from [\DDR] (or from [\Survey]) that every subequation $F\ss J^2X$ on a manifold $X$ satisfies positivity and 
negativity conditions.  In local coordinates $x=(x_1, ... , x_n)$ where we have  a canonical trivialization
$J_x \rn = \bbr\oplus \rn\oplus \Symn$, this means that  $F$ and $\Int F$ are preserved under translations by 
$(r, 0, P)$ with $r\leq 0$ and $P\geq 0$.


\Def{\DD.12. (The Coordinate Ball Condition)}   We say  that a reduced convex cone subequation $F$  satisfies 
{\bf the coordinate ball condition} 
if  each point $x_0\in X$
has a local coordinate neighborhood $U$ and $c>0$ such that the reduced jet
$$
(0, c I)_{x_0} \in \Int \left\{F\bigr|_U\right\}.
$$

Note that if  this holds for one coordinate neighborhood about a point $x_0$, then it will hold for 
all coordinate neighborhood about $x_0$. To see this,   note that after a change of coordinates fixing
$x_0$ the reduced jet $(0, c I)_{x_0}$ becomes $(0, P)$ for $P>0$. Now there exists a $P_0\geq0$
such that $P+P_0 = c' I$ for $c'\geq c$.  Since $P\in \Int F$, we have $P+P_0\in \Int F$.

\Lemma{\DD.13} {\sl
Suppose that $F$ is a reduced convex cone subequation on a manifold $X$.
Let $U$ be a local coordinate system about a point $x_0\in X$.
If $F$ satisfies the coordinate ball condition, then there exists a $c_0>0$
 such that  for each $c\geq c_0$, the  $\ci$ defining function
$$
\rho(x) \ \equiv\  \smfrac c 2 \left( |x-x_0|^2 - \d^2\right)
$$
for the boundary of the coordinate ball $B_\d(x_0) \equiv \{x : |x-x_0|<\d\}$,
is strictly $F$-subharmonic on $B_\d(x_0)$ for all  $\d>0$  sufficiently small.
}

\noindent
{\bf Proof.}  By the note after Definition \DD.12 our condition holds in
any  given local coordinates.   The reduced 2-jet of $\rho$ at $x$ equals $J^{2, {\rm red}}_x \rho = c (x-x_0, I)$.
At $x=x_0$ this jet belongs to $(\Int F)_{x_0}$ by the hypothesis.  
Hence, it belongs to $\Int F$  for $x$ near $x_0$, proving that $\rho$ is $F$-strict on 
$B_\d(x_0)$ for $\d>0$ small.\qed

The upshot of this observation is that in applying any of the above methods for verifying the 
local hypothesis, we need only consider domains $\O$ which are very small balls in some local
coordinate system.

\noindent
{\bf Proof of Theorem \BB.5.}  In this theorem we are assuming that there is a covering family of
local coordinates on $X$ such that for all sufficiently small balls $\O$ in these coordinates,
the $\ci$ homogeneous $F$-Dirichlet problem is uniquely solvable on $\O$, that is, Definition \DD.6 holds.
We are also supposing that $F$ satisfies the coordinate ball condition in Definition \DD.12.
By  the note following Lemma \DD.13, we know that Lemma \DD.13 applies to all sufficiently small
balls in this covering family of coordinates.  That is, each such ball has a smooth strictly $F$-subharmonic
defining function $\rho$.

Now we are able to apply Proposition \DD.7 to conclude that approximate $\ci$-regularity for the $F$-(DP)
holds on all sufficiently small balls $\O$ about each point in our family of coordinates.
By  Theorem \BB.2$'$ (located after proof of Lemma \DD.5),
 we conclude that  for all $u\in  F_{\rm cont}^{\rm strict}(X)$,
Richberg approximation holds for $u$ on $X$.\qed

\noindent
{\bf Proof of Theorem \BB.6.}  The argument here is parallel to the one above.
In this theorem we are assuming that there is a covering family of
local coordinates on $X$ such that for all sufficiently small balls $\O$ in these coordinates,
the $\ci$ inhomogeneous $F$-Dirichlet problem  (IHDP) is uniquely solvable on $\O$, that is, Definition \DD.8 holds.
We are also supposing that the monotonicity subequation $M$ 
satisfies the coordinate ball condition in Definition \DD.12.
By  the note following Lemma \DD.13, we know that Lemma \DD.13 applies to all sufficiently small
balls in this covering family of coordinates.  That is, each such ball has a smooth strictly $M$-subharmonic
defining function $\rho$.

Now we are able to apply Proposition \DD.9 to conclude that approximate $\ci$-regularity for the $F$-(DP)
holds on all sufficiently small balls $\O$ in our family of coordinates.
We then apply   Theorem \BB.2$'$ to complete the proof as above.\qed

\Remark {\DD.14. (A Special Case of Theorem \BB.6)}
Consider the  case  where
 $F$ is a reduced convex cone subequation with compatible operator $f\in C(F)$. 
 In this case we can take $M=F$ and the hypothesis (1) before Definition \DD.8 and be
 replaced by the super additive condition:
$$
\text{
For all $J, J' \in F$, one has that $f(J+J')\ \geq\ f(J) + f(J')$
}
\eqno{(1)'}
$$
because this implies (1) with $\d=1$ and $g=f$ (since $f(J)\geq0 \ \forall J\in F$).
Note also  that (1) $\Rightarrow$ (1)$'$ with   $\d = 2$ and $g=f$.
This is because  $f(J)+f(J') \geq 2f(J)$ and $f(J)+f(J') \geq 2f(J')$, when added up,
 give (1)$'$.

\vskip.3in

\noindent
{\headfont \EE. Smoothing the Maximum Function.}

In this section we discuss the proof of Properties \CC.1 and Lemma \CC.2.
The maximum function $M(t) \equiv \max\{t_1,...,t_m\}$ has the obvious properties
that $M$ is convex, invariant under permutations and, with $e=(1,...,1)$, 
$$
M(t+se) \ =\ M(t) + s.
\eqno{(\EE.1)}
$$
These same properties carry over to $M_\e(t)$. This gives Property (1).  For Property (2)
note that convolution preserves (\EE.1) since $s * \vf  = s$ by (\CC.2).  
Taking ${d\over ds}\bigr|_{s=1}$ of both sides of 
the equation $M_\e(t+se) \ =\ M_\e(t) + s$ yields $\sum {\partial M_\e \over \partial t_j}=1$.
Properties (3) and  (4) are equally straightforward. This proves Properties \CC.1.

\noindent
{\bf Proof of Lemma \CC.2.}  The following is a straightforward calculation.

\Lemma{\EE.1} {\sl
With $u_1,...,u_m\in C^\infty(X)$ arbitrary, and $w\equiv M_\e(u_1,...,u_m)$, the 2-jet of $w$
equals
$$
J^2 w \ =\ \sum_{j=1}^m {\partial M_\e \over \partial t_j} J^2 u_j  +  (E_\e, 0, P)
\eqno{(\EE.2)}
$$
where
$$
E_\e \ =\  w - \sum_{j=1}^m {\partial M_\e \over \partial t_j}  u_j  
\and
P\ \equiv \  \sum_{i,j=1}^m {\partial^2 M_\e \over \partial t_i \partial t_j}  D u_i  \circ D u_j.
\eqno{(\EE.3)}
$$
}
\smallskip

Since $M_\e$ is convex, $P\geq0$.  Now the fibres of $F$ and of $F^{\rm strict}_g$
at a point $x\in X$ are convex.  Hence, $J_x^2 u_j\in F_x^{g(x)}$ for all $j=1,...,m$
implies that the convex combination 
$\sum_j{\partial M_\e \over \partial t_j} J^2 u_j  \in F_x^{g(x)}$.  By positivity, 
$\sum_j{\partial M_\e \over \partial t_j} J^2 u_j +(0,0, P)  \in F_x^{g(x)}$.
The error $E_\e \equiv M_\e(t) - \sum_j t_j {\partial M_\e \over \partial t_j}$
satisfies
$$
-M(\e) \ \leq\ E_\e(t) \ \leq \ M(\e).
\eqno{(\EE.4)}
$$
This proves that $J_x^2 w\in F_x^{g(x)-M(\e)}$.

The estimate (\EE.4) is verified as follows.
First compute that 
$$
\begin{aligned}
E_\e(t) \  &=  \ {d\over dr} \biggr|_{r=1} \left(rM_\e(t) - M_\e(rt)   \right) \\
 &=  \  {d\over dr} \biggr|_{r=1}  \int_{\bbr^m} \biggl(r M(t+\e y) - M(rt+ey)\biggr) \, \vf(y) \, dy \\
&=  \  {d\over dr} \biggr|_{r=1}  \int_{\bbr^m} \biggl(M(rt+r\e y) - M(rt+ey)\biggr) \, \vf(y) \, dy \\
&=  \  \int_{\bbr^m} \sum_{j=1}^m \e_j  y_j {\partial M\over \partial t_j} (t+\e y) \, \vf(y)\, dy.
\end{aligned}
$$
Since $y_j \in \supp (\vf) \Rightarrow |y_j|\leq 1$, we have $-M(\e) \leq \e_j y_j \leq M(\e)$.
Therefore, the convex combination $\sum \e_j y_j {\partial M\over \partial t_j} $ in the
integral lies between $-M(\e)$ and $M(\e)$.\qed

\medskip
\noindent{\headfont \FF.\  An Example -- Non-convex Subequations.}

The problem of smoothing $F$-subharmonic functions is also  interesting for subequations which are not necessarily convex.  Diederich and Forn\ae{}ss [\DF]  show that it is impossible,  for  general $F$, to regularize the maximum of a finite number of smooth $F$-subharmonic functions. On the other hand, they prove that it is possible for the case of smooth $n$-convex functions, that is, for the subequation $F_n=\widetilde F_1=\widetilde{PSH}$, which is not convex.\footnote{\ Here $n$ is the complex dimension of the manifold and the subequation $F_n$ is defined by requiring that at least
one eigenvalue of the complex hessian is $\geq0$.} The following Proposition shows that convexity of a subequation is not a necessary condition for Richberg approximation.   

\Prop{\FF.1}  {\sl
Let $X$ be a complex manifold.
For any continuous strictly $n$-convex function $u$ on $X$,
and any $h\in C(X), \ h>0$,
 there exists $$w\in C^\infty(X)\cap F_n^{\rm strict}(X)$$ which satisfies:
$$
u\ <\ w\ <\ u+h 
$$
on $X$.
}

{\bf Proof.} We show first that (locally) approximate $\ci$ regularity for the (homogeneous)
(DP) holds, that is,
any $F_n$-harmonic function is the limit of a sequence of smooth strictly $n$-convex functions. Let $\Omega\subset \mathbb{C}^n$ be strictly pseudoconvex and let $H\in C(\Omega)$ be  $F_n$-harmonic. Then $-H$ is a plurisubharmonic function and by [\CKNS] there is the sequence $w_k\subset\Omega$ of  smooth plurisubharmonic functions which satisfies the following conditions:\\

(i) \ \ $w_k$ converge uniformly to $-H$, \ \  and\ \ 

(ii)\ \  $\det\left(\frac{\partial^2w_k}{\partial z_p\partial\bar z_q}\right)_{p,q=1}^n=\left(\frac{1}{2k}\right)^n$.

\noindent
In particular there is an eigenvalue of $\left(\frac{\partial^2w_k}{\partial z_p\partial\bar z_q}\right)_{p,q=1}^n$ which is smaller than $\frac{1}{k}$ and therefore the function $w_k-\frac{1}{k}|z|^2$ is nowhere plurisubharmonic. Thus the sequence $v_k=\frac{1}{k}|z|^2-w_k$ is a sequence of smooth strictly $n$-convex functions which converge uniformly to $H$.

By Lemma \DD.5, approximate $\ci$ regularity for the  (DP) implies quasi $\ci$ approximation.
(See [\DDI] for the step in the proof --  Lemma 3.1 -- where solving the (DP) is required.)

 Finally, we can choose $\Omega_\alpha$, $K_\alpha$ and $v_\alpha$ as in the proof of Theorem \BB.2. The function 
 $$
 \tilde w(z)=\sup\{v_\alpha(z):z\in\Omega_\alpha\}>u
 $$
  is, in  the terminology of [\DF], $n$-convex with corners, and by the Diederich-Forn\ae{}ss approximation result (Theorem 1 in [\DF]) it can be approximated by $w$ satisfying the statement in Proposition \FF.1.\qed

 \vskip.3in

\def\item{}
\vskip .3in

\centerline{\headfont References}

\noindent
\item{[\CKNS]}     L. Caffarelli, J. J. Kohn, L. Nirenberg and  J. Spruck,
  {\sl The Dirichlet problem for non-linear second
order elliptic equations II: Complex Monge-Amp\`{e}re, and
uniformly elliptic equations}, Comm.\ Pure Appl.\ Math.\  {\bf 38} (1985),
209-252,

\smallskip

\noindent

\noindent
\item{[\CTW]} J. Chu, V. Tosatti and B. Weinkove,  {\sl
The Monge-Ampre equation for non-integrable almost complex structures}, 
J. Eur. Math. Soc. {\bf 21} (2019), no.7, 1949-1984. ArXiv:1603.00706.

\smallskip

\noindent
\item{[\D]} J.-P.  Demailly, {\sl Complex analytic and differential geometry}.  An e-book available at:
http://www-fourier.ujf-grenoble.fr/~demailly/documents.html.

\smallskip

\noindent
\item {[\DF]}  K. Diederich, I. E. Forn\ae{}ss,  {\sl Smoothing q-convex functions and vanishing theorems}, Invent. Math. 82 (1985), no.\ 2, 291-305.

\smallskip
\noindent
\item{[\GW]} R. E. Greene and H. Wu, {\sl $\mathcal{C}^\infty$-approximations of convex, subharmonic, and plurisubharmonic functions},
Ann.\ Scient.\ Ec.\ Norm.\ Sup. 12, 47–84 (1979).

\smallskip
\noindent
\item{[\GN]} D. Gu and N. C. Nguyen {\sl The Dirichlet problem for a complex Hessian equation on compact Hermitian manifolds with boundary}, Ann. Sc. Norm. Super. Pisa Cl. Sci. (5) Vol. XVIII (2018), 1189-1248.

\smallskip

 \noindent
\item {[\PUP]}     F. R. Harvey and H. B. Lawson, Jr.  ,  {\sl  Plurisubharmonicity in a general geometric context},
  Geometry and Analysis {\bf 1} (2010), 363-401. ArXiv:0804.1316

\smallskip

 \noindent
\item {[\DDI]}   \ \----------,  {\sl  Dirichlet duality and the non-linear Dirichlet problem},    Comm. on Pure and Applied Math. {\bf 62} (2009), 396-443.

\smallskip

 \noindent
\item{[\DDR]}  \ \----------, {\sl Dirichlet Duality and the Nonlinear Dirichlet Problem on Riemannian Manifolds},  J. Diff. Geom. {\bf 88} (2011), 395-482.   ArXiv:0912.5220.

 \noindent
\item {[\GPSH]} \ \----------,  {\sl  Geometric plurisubharmonicity and convexity - an introduction},
  Advances in Math.  {\bf 230} (2012), 2428-2456.   ArXiv:1111.3875.

 \smallskip
 \noindent
\item {[\Survey]}  \ \----------,   {\sl  Existence, uniqueness and removable singularities
for nonlinear partial differential equations in geometry},\ 
 pp. 102-156 in ``Surveys in Differential Geometry 2013'', vol. 18,  
H.-D. Cao and S.-T. Yau eds., International Press, Somerville, MA, 2013.
ArXiv:1303.1117.

\smallskip

\noindent
\item {[\Bellman]}  \ \----------,  {\sl The equivalence of  viscosity and distributional
subsolutions for convex subequations -- the strong Bellman principle}, Bulletin Brazilian Math.\ Soc.\ {\bf 44} No.\ 4 (2013),  621-652.  ArXiv:1301.4914.

\smallskip

 \noindent
\item {[\Lagrangian]}  \ \----------, {\sl  Lagrangian potential theory and a Lagrangian equation of  Monge-Amp\`ere type}, pp. 217- 257 in Surveys in Differential Geometry, 
 edited by H.-D. Cao, J. Li, R.  Schoen and S.-T. Yau, {\bf  22},  International Press, Somerville, MA, 2018.
  ArXiv:1712.03525.

\smallskip



 \noindent
\item {[\HLP]}   F. R. Harvey, H. B. Lawson, Jr.  and S. Pli\`s,  {\sl  Smooth approximation of plurisubharmonic functions on almost complex manifolds},  Math. Ann.  {\bf 366}, issue 3, (2016), 929-940.     ArXiv:1411.7137.

\smallskip

 \noindent
\item{[\Plis]}  \ \  S. Pli\'s,  {\sl The Monge-Amp\`ere equation on almost complex manifolds},  Math. Z.
{\bf 276} (2014), no. 3-4, 969-983.
 \smallskip

\noindent
\item {[\Pl]}  \ \----------,   {\sl On regularization of
$J$-plurisubharmonic functions}, C.R. Acad. Sci. Paris, Ser.I (2014),
http://dx.doi.org/10.1016/j.crma.2014.11.001. 

\smallskip

 \noindent
\item{[\Pll]}   \ \----------,  {\sl  Monge-Amp\`ere operator on four dimensional  almost complex manifolds}, 
 J. Geom. Anal. {\bf 26}  (2016), no. 4, 2503Ð2518.  
ArXiv: 1305.3461.

\smallskip
 \noindent
\item{[\Plll]}  \ \----------,   {\sl The smoothing of $m$-subharmonic functions},
 arXiv:1312.1906.


\smallskip
 \noindent
\item{[\Ri]}
 R. Richberg, {\sl Stetige streng pseudokonvexe Funktionen}, Math. Ann. {\bf 175}
(1968), 257-286.


\smallskip
\noindent
\item{[\V]} M. Verbitsky, 
{\sl Plurisubharmonic functions in calibrated geometry and q-convexity}, 
Math.\ Z.\ 264 (2010), no.\ 4, 939-957.

 \end{document}